\documentclass{amsart}
\usepackage{amssymb,amsmath,amsthm,latexsym,amsfonts,epsfig,euscript}
\newtheorem{fed}{Definition}[section]
\newtheorem{teo}[fed]{Theorem}
\newtheorem{lem}[fed]{Lemma}
\newtheorem{cor}[fed]{Corollary}
\newtheorem{pro}[fed]{Proposition}
\newtheorem{rem}[fed]{Remark}

\newtheorem{num}[fed]{}

\def\bdem{\begin{proof}}
\def\edem{\end{proof}}

\newfont{\bb}{msbm10}

\def\orto{^\perp}
\def\inc{\subseteq}

\def\inv{^{-1}}
\def\*A{\#\sb A}

\def\H{{\mathcal H}}
\def\iiff{if and only if }

\def\ca{L(\H ) }
\def\cam{L(\H )^+ }

\def\cH{{\mathcal H}}

\def\cQ{{\mathcal Q}}
\def\cR{{\mathcal R}}

\def\cS{{\mathcal S}}
\def\cT{{\mathcal T}}
\def\cM{{\mathcal M}}
\def\cN{{\mathcal N}}
\def\cW{{\mathcal W}}
\def\cK{{\mathcal K}}
\def\cD{{\mathcal D}}
\def\cG{{\mathcal G}}

\def\PA{P_{A, \cS}}

\def\QA{Q_{A,\cS} }

\def\PAS{\mathcal
 P(A, \cS)}

\def\rai{^{1/2}}

\def\api{\langle}
\def\cpi{\rangle}

\def\noi{\noindent}

\def\bm{\left(\begin{array}}
\def\em{\end{array}\right)}
\def\ben{\begin{enumerate}}
\def\een{\end{enumerate}}
\def\beq{\begin {equation}}
\def\endeq{\end{equation}}
\def\barr{\begin{array}}
\def\earr{\end{array}}

\def\inv{^{-1}}

\def\H{{\mathcal H}}
\def\lh{{L(\H)}}
\def\lh+{{\lh^+}}

\def\com{$(A, \cS)$ is compatible}

\def\ov{\overline}
\def\sh{\sharp}
\def\wtil{\widetilde}

\begin{document}
%\fontsize {14}{17}\selectfont
\title{Projections in Operator ranges 
\footnote{2000 AMS-MSC: Primary 46C07, 47A62, 46C05.
}
}

\author {G. Corach 
 A. Maestripieri and D. Stojanoff 
}

\thanks{Partially supported by CONICET (PIP 2083/00), UBACYT I030 and ANPCYT (PICT03-9521)}

\maketitle

\vskip.5truecm
\begin{abstract}{
If $\H$ is a Hilbert space, $A$ is a positive bounded linear operator on $\cH$ and
$\cS$ is a closed subspace of $\cH$, the relative position between
$\cS$ and $A^{-1}(\cS \orto)$ establishes a notion of compatibility. 
We show that the compatibility of $(A,\cS)$ is equivalent to the existence 
of a convenient orthogonal projection
in the operator range $R(A^{1/2})$ with its 
canonical Hilbertian structure. 
}\end{abstract}

\newdimen\normalbaselineskip
%\normalbaselineskip=16pt
\normalbaselines \vskip 1truecm
%\fontsize {14}{17}\selectfont

\section{Introduction}
Oblique projections are becoming an important tool in several
areas of mathematics, statistics and engineering. This phenomenon 
is illustrated in many papers on integral equations, iterative methods
in numerical linear algebra,
signal processing, linear regression, just to mention a sample; in 
\cite{[CMS4]} the reader can find an extensive list of papers
on these applications. In a recent 
series of papers \cite{[Sz]},
\cite{[CMS2]}, \cite{[CMS3]}, \cite{[CMS4]} the set of oblique
projections is studied according to different inner and semi-inner
products which orthogonalize them. This is the way in which a
certain notion of compatibility arises. A positive (Hermitian semidefinite)
operator $A$ on
a Hilbert space $\cH$ and a closed subspace $\cS$ of $\cH$ are
said to be {\bf compatible} if there exists a projection $Q$ in
$\cH$ with range $\cS$ such that $AQ=Q^*A$. This equality means
that $(Qx , y )_A=(x , Qy ) _A$ $\forall x ,y \in\cH$ if $( u , v
)_A:=\langle A u , v \rangle$ where $ u , v \in\cH$ and $\langle
~,~ \rangle$ denotes the inner product on $\cH$. Observe that $(
~,~ )_A$ is, in general, a semi-inner product, because $A$ is allowed
to have a non trivial nullspace. If the pair $(A,\cS)$ is
compatible then a distinguished element $P_{A,\cS}$ in
$$
P(A,\cS)=\{Q\in L(\cH):Q^2=Q, ~Q\cH=\cS, ~AQ=Q^*A\}
$$
can be defined with certain optimal properties.

On the other hand, given Hilbert spaces $\cH, \cK$ the range of a
bounded linear operator $T: \cH \to \cK$ can be naturally
given a Hilbert space structure, by means of the inner product
$\langle Tx ,Ty \rangle_T=\langle x_1 ,y_1 \rangle$, $x ,y \in \cH$, 
where $x_1$ (resp. $y_1$) denotes the orthogonal projection of $x$ (resp. $y$) 
to the closure of $R(T)$ in $\cK$. 
These Hilbert spaces $B(T) = (R(T),\langle~,~\rangle_T)$
play a significant role in many areas, in particular in the de
Branges complementation theory. The reader is referred to the
books by de Branges and Rovnyak \cite{[Br]} and Ando \cite{[An]}
for systematic expositions of this theory. The main goal of this paper
is to determine the compatibility of a pair $(A,\cS)$ by checking the
existence of a convenient orthogonal projection in the space
$B(A^{1/2})=(R(A^{1/2}),\langle ~,~ \rangle_{A^{1/2}})$. This
approach allows us to see the oblique projection $P_{A,\cS}$ as a
true orthogonal projection (acting, of course, on a different
Hilbert space, namely $B(A^{1/2})$). Let us describe more
precisely these concepts and results. 
 Section 2 collects
some notations and a description of a theorem by R. G. Douglas which is one of
the main tools of this paper. Douglas theorem studies the existence and
uniqueness
of solutions of operator equations like $AX = B$, for  operators
$A, B$ between Hilbert spaces. 
Section 3 starts with
a survey of known results on compatibility and on the form of a
distinguished projection $P_{A,\cS}$ with the properties mentioned
above. Some proofs of these results can be found in \cite{[Sz]}, 
\cite{[CMS3]} and \cite{[CMS4]}.
In addition, we present new characterizations of compatibility;
some of them are quite technical but they will be needed later, in the
sections dealing with operator ranges. 
Section 4 contains a description of the Hilbertian structure on an
operator range. Here, the references are a paper by Fillmore and
Williams \cite {[FW]} and the books by Ando \cite{[An]} and de Branges
and Rovnyak \cite {[Br]}, in
addition to a paper by Dixmier \cite{[Di]}. The particular
operator range we are interested in is $R(A^{1/2})$, i.e., the
range of the positive square root of a fixed positive operator
$A$. Moreover, we need to characterize the closure and the
orthogonal complement of a subspace in $B(A^{1/2})$ and the
algebra of all bounded operators on $\cH$ which can be extended,
after a convenient reduction modulo the nullspace $N(A)$ of $A$,
to $B(A^{1/2})$. In this section we slightly extend some results
by Barnes \cite{[Bar]} who studied the case of an injective
operator $A$; however, Barnes' goal is different from ours, namely,
he studies the spectral
properties of an operator when it is seen in $B(A)$ or in $B(A^{1/2})$. 
Finally, Section 5 contains
a characterization of
the compatibility of a pair $(A,\cS)$ in terms of
certain decompositions of $B(A^{1/2})$. Moreover, it is proven
that if $(A,\cS)$ is compatible then the distinguished projection
$P_{A,\cS}$ can be extended (in the sense mentioned above) to
$B(A^{1/2})$, and conversely. Also, it is shown that the orthogonal
projection $P_{\cW}$ onto a closed subspace $\cW$ of 
$B(A^{1/2})$ comes from an operator on $\cH$ if
and only if $(A,\cS)$ is compatible, where $\cS$ is a closed subspace of
$\cH$ such that $A(\cS)$ is dense in $\cW$ (in the topology of $B(A^{1/2})$).

\section{Preliminaries}
In what follows $\cH$ and $\cK$ denote  Hilbert spaces, 
%with inner products $\langle ~ , ~ \rangle$, 
$L(\cH, \cK)$ is the Banach space of bounded linear
operators from $\cH$ to $\cK$,  $L(\cH)$ is the algebra $L(\cH, \cH)$
and $\cam$ denotes the cone of positive operators on $\cH$.  The set of all (bounded linear)
projections in a Hilbert space is denoted by $\cQ$.
For
any $W \in \ca$, the range and the nullspace of $W$ are respectively denoted 
by $R(W)$ and $N(W)$.  Given a closed subspace $\cM$ of $\cH$, $P_{\cM}$ denotes 
the orthogonal projection onto $\cM$. If  
$W\in L(\cH, \cK)$ has closed range, then
the Moore-Penrose pseudoinverse of $W$, denoted by $W^\dagger$, belongs to
$L(\cK, \cH)$ and it is characterized by the properties $WW^\dagger=P_{R(\cW)}$
and $W^\dagger W=P_{R(\cW^*)}$ (see \cite{[DW]},  \cite{[BG]} and \cite {[De1]}
for more properties and applications of $W^\dagger$).

We state the theorem by R. G. Douglas \cite{[Do]}, \cite{[FW]} mentioned
in the introduction which will be used in several parts of the paper.

\medskip
\begin{teo} Given
Hilbert spaces $\cH$, $\cK$, $\cG$ and operators $A\in
L(\cH,\cG)$, $B\in L(\cK,\cG)$ then the following conditions are
equivalent:

i) the equation $AX=B$ has a solution in $L(\cK,\cH)$;

ii) $R(B)\subseteq R(A)$;

iii) there exists $\lambda >0$ such that $BB^*\le \lambda AA^*$.
In this case, there exists a unique $D\in L(\cK,\cH)$ such that
$AD=B$,  $R(D)\subseteq \ov{R(A^*)}$, and $N(D) = N(B)$; moreover,
$\|D\|^2=\inf\{\lambda >0:~ BB^*\le \lambda AA^*\}$.
We shall call $D$ the {\bf reduced solution} of $AX=B$.
\end{teo}

As a consequence of Douglas theorem and the properties of the Moore-Penrose
pseudoinverses, it follows that if $R(A)$ is closed and $R(B)\subseteq R(A)$ then
$A^\dagger B$ is the reduced solution of $AX=B$.

\section{Oblique projections}

Given $A\in L(\cH)$, the functional
$$
( \ , \  )_A : \H \times \H \to \mathbb C \ , \quad ( x, y
)_A = \api Ax, y \cpi , \quad x, y \in \H
$$
is an equivalent inner product on $\H$ if and only if $A$ is a
positive invertible operator on $\cH$. If $A \in \cam$, then
$( \ , \
)_A $ is a Hermitian sesquilinear form which is positive semidefinite,
i.e., a semi-inner product on $\cH$.
For a subspace $\cM$ of $\H$ it is easy to see that
$$
\{ x \in \H : ( x, y )_A = 0 \ \forall y \in \cM \} = (A
\cM ) \orto = A\inv (\cM \orto ).
$$
Given $W \in \ca$, an $A$-adjoint of $W$ is any $V \in \ca $ such
that $ ( Wx, y )_A = ( x, V y )_A  , \quad \hbox{for
all} \ x, y \in \H , $ i.e., $AW = V^*A$. We are interested
in projections $Q \in \cQ$ which are $A$-Hermitian, in the sense
that $AQ=Q^*A$. 

From now on, we fix $A \in \cam$ and a closed
subspace $\cS$ of $\cH$ and abbreviate $P = P_{\cS}$. The first result is due 
to M. G. Krein
\cite{[K]}. There is a recent proof of it in \cite{[Sz]}.
\medskip
\begin{lem}[Krein] \label{L2} Let $Q$ be a projection with $R(Q)=\cS$.
Then $Q$ is $A$-Hermitian
if and only if $N(Q)\inc A\inv (\cS\orto)$.
\end{lem}

\medskip
Recall from the introduction the set $P(A,\cS)$ of all $A$-Hermitian 
projections with
fixed range $\cS$, i.e. $P(A,\cS)=\{Q\in\cQ:R(Q)=\cS$ and
$AQ=Q^*A\}$. The pair $(A,\cS)$ is {\bf compatible} if the set
$P(A,\cS)$ is non empty.

Observe that it follows from lemma \ref{L2} that if a projection
$Q$ has range $\cS$ then $Q\in \PAS$ \iiff $N(Q) \inc A\inv
(\cS\orto )$, so that \com \ \iiff $\H = \cS + A\inv (\cS\orto )$.
In this case, $\cH=\cS\oplus(A^{-1}(\cS\orto)\ominus\cN)$, where
$\cN=\cS\cap A^{-1}(\cS\orto)=\cS\cap N(A)$ and there exists a
unique projection $P_{A,\cS}$ with range $\cS$ and nullspace
$A^{-1}(\cS\orto)\ominus\cN$. It is elementary to check that
$P_{A,\cS}\in P(A,\cS)$. At the end of section 3 we shall mention some
optimal properties of $P_{A,\cS}$.
\medskip

\begin{rem}\rm
In \cite{[Ba]}, Baksalary and Kala studied, in the matrix case,
the existence of $P_{A,\cS}$
under the additional hypothesis of the invertibility
of $A$. In \cite{[HN]}, Hassi and Nordstr\"om determined conditions on a
Hermitian not necessarily invertible operator $A$, under which
the set $P(A,\cS)$ is a singleton. They also proved some
least-square-type results for indefinite inner products.
In \cite{[PW]}, Z. Pasternak-Winiarski studied, for $A$ invertible,
the analiticity of the map $A \to P_{A, \cS}$ (see also \cite {[ACS1]}
for shorter proofs and related results).
\end{rem}

\bigskip
Consider  the matrix representation of $A$ in terms of the
orthogonal projection $P_{\cS}$ onto $\cS$, namely,

\begin{equation}
A=\begin{pmatrix} a & b\\ b^* &c\end{pmatrix};
\end{equation}
this means that $a\in L(\cS)$, $b\in L(\cS\orto,\cS)$, $a\in
L(\cS\orto)$ and $Ax=as + bs \orto + b^*s+cs\orto$ if $x=s+s\orto$
is the decomposition of $x\in\cH=\cS\oplus\cS\orto$. If $Q\in\cQ$
and $R(Q)=\cS$ then there exists $X \in L(\cS\orto,\cS)$ such that 
the matrix representation of $Q$ in terms of 
$P_{\cS}$ is $Q=\begin{pmatrix} 1 &X\\ 0 & 0\end{pmatrix}$. 
It is easy to see that the condition $AQ=Q^*A$ is equivalent
to the equation $aX=b$. Then $(A,\cS)$ is compatible if and only
if the equation $aX=b$ admits a solution. Applying Douglas
theorem, this is equivalent to $R(b)\subseteq R(a)$
(or $R(PA)\subseteq R(PAP)$). Consider the reduced solution $d$ of
$aX=b$. It easily follows that $P_{A,\cS}=\begin{pmatrix} 1 & d
\\ 0 & 0 \end{pmatrix}$ (see \cite{[Sz]},
for a proof of these facts). Observe that, if $A$ is invertible, then
$P_{A,\cS}= P\Big(PAP + (I-P)A(I-P)\Big)^{-1}A$. For many results in
the case of invertible $A$,
the reader is referred to \cite {[PW]}, \cite {[ACS1]} and \cite {[Sz]}.

Some basic conditions for the compatibility of the pair $(A,\cS)$
can be found in \cite{[Sz]}, \cite{[CMS3]}, \cite{[CMS4]} as well
as formulas for the elements of $P(A,\cS)$, if $(A,\cS)$ is
compatible.

In what follows we give new characterizations of compatibility;
also, we express the distinguished element $P_{A,\cS}$ of
$P(A,\cS)$ as the solution of certain Douglas-type equations.

\medskip
\noi
If $\cS\cap N(A)=\{0\}$ then the compatibility of $(A,\cS)$ can be
easily checked. In fact:
\medskip
\begin{pro}\label{5}
Consider $A\in L(\cH)^+$ such that $\cS\cap N(A)=\{0\}$. Then the
following conditions are equivalent:

i) $\cH=\cS \oplus A(\cS)^\bot$, i.e., $(A,\cS)$ is compatible;

ii) $\overline{A(\cS)}\oplus \cS^\bot$ is closed;

iii) $\overline{A(\cS)}\oplus \cS^\bot=\cH$.
\end{pro}

\bdem \ \ben

\item[i)]$\rightarrow$ ii) We use the general fact that if
$\cM,\cN$ are closed subspaces, $\cM +\cN$ is closed if and only
if $\cM^\bot+\cN^\bot$ is closed (see theorem 4.8 of \cite{[Ka]}):
if $\cS\oplus A(\cS)^\bot=\cH$, a fortiori $\cS+A(\cS)^\bot$ is
closed. Then $\cS^\bot+\overline {A(\cS)}$ is closed. Besides
$\cS^\bot\cap \overline{A(\cS)}=(\cS+A(\cS)^\bot)^\bot=\{0\}$.

\item[ii)]$\rightarrow$ iii) If $\cS^\bot+\overline{A(\cS)}$ is
closed, $\cS^\bot+\overline{A(\cS)}
=\overline{\cS^\bot+\overline{A(\cS)}}=(\cS\cap
A(\cS)^\bot)^\bot=(\cS\cap N(A))^\bot=\cH$.

\item[iii)] $\rightarrow$ i) is similar. \een \edem

The closure condition of  part ii) is equivalent to an angle
condition. In fact, the sum of two closed subspaces is
closed if the angle determined by them is non zero. The reader is
referred to \cite{[Ka]}, \cite{[BG]}, \cite{[De]}, \cite {[HJ]}
for  nice surveys on angles in
Hilbert spaces,  and to \cite{[Sz]}, \cite{[CMS4]} for particular
details concerning compatibility. In particular, in \cite{[CMS4]}
it is proven that $(A, \cS )$ is compatible if and only if the angle between 
$\cS\orto$ and the closure of $A(\cS)$ is non zero.
\medskip

\noi
Next, we state  a chain of necessary conditions for compatibility.

\medskip

\begin{pro} \label{R3.4}
Consider the following conditions:
\ben 
\item The pair $(A,\cS)$ is compatible. 
\item $A(\cS)$ is
closed in $R(A)$. \item $A^{-1}(\overline{A(\cS)})=\cS +N(A)$.
\item  $A^{1/2}(\cS)$ is closed in $R(A^{1/2})$. \item $\cS+N(A)$
is closed. \item $P_{\overline{R(A)}}(\cS)$ is closed. 
\een 
Then  1 $\rightarrow$ 2 $\rightarrow$ 4
$\rightarrow$  5,  2 $\leftrightarrow$ 3 and 5 $\leftrightarrow$
6.
\end{pro}

\bdem \ 

\ben 

\item [1 $\rightarrow$ 2:] Observe that  \com\ if and only if
$R(A)=A(\cS)\ +\ (\cS^\perp\cap R(A))$. Consider $z\in
\overline{A(\cS)}\cap R(A)$; then  there exists a sequence
$\{s_n\}$ in $\cS$ such that $As_n\rightarrow z$ and there exist $s\in \cS$ and
$y \in \cH$ such that $Ay\in \cS^\perp$  and $z=As+Ay$. Since $\api As_n, w\cpi \ =
0$ for every $w\in A\inv (\cS ^\perp )$, then $\api z, w\cpi =0$
for every $w\in A\inv (\cS ^\perp )$. Thus, $0=\api z,y\cpi= \api
s,Ay\cpi +\api Ay ,y\cpi= \api Ay,y\cpi= \|A^{1/2}y\|^2$ and $y\in
N(A)$. Therefore, $z=As\in A(\cS)$.

 \item[2 $\rightarrow$ 4:] Consider $z\in\overline{ A^{1/2}(\cS)}\cap R(A^{1/2})$;
 then $z=A^{1/2}x$
 for some $x\in\H$ and there exists a sequence $\{s_n\}$ in $ \cS$ such
 that $A^{1/2}s_n\rightarrow A^{1/2}x$ ;
 then $As_n\rightarrow Ax$ so that
 $Ax=A^{1/2}z\in\overline{A(\cS)}\cap R(A)=A(\cS)$.
 Then $z\in (A^{1/2}(\cS)+N(A))\cap R(A^{1/2})$ so that $z\in A^{1/2}(\cS)$.

\item[2 $\leftrightarrow$ 3:] Observe that $A(\cS)$ is closed in $R(A)$ if and only if
$\overline{A(\cS)}\cap R(A)=A(\cS)$; but this is equivalent to
$A^{-1}(\overline{A(\cS)})=A^{-1}(A(\cS))$ , i.e.,
$A^{-1}(\overline{A(\cS)})= \cS+N(A)$.

\item[4 $\rightarrow$  5:] It is easy to prove that if 
$A^{1/2}(\cS)$ is closed in $R(A^{1/2})$, then 	 $A^{-1/2} (\overline{A^{1/2}(\cS)})$ $=\cS+N(A)$. 
So that $\cS+N(A)$ is closed.

 \item[5 $\leftrightarrow$ 6:] It is a general result that if $\cM$ and $\cN$ are
 closed subspaces  then $\cM +\cN$ is closed if and only if $P_{\cN ^\perp}(\cM)$ is closed
 (see \cite{[Ka]} or \cite{[De]}.
\een
\edem

\noi
In \cite {[Sz]} it is shown that all  conditions  above are equivalent
if $R(A)$ is closed.
The next technical result will be used in the following sections.

\medskip
\begin{cor}\label{R3.6}If $\cS + N(A)$ is closed then the following conditions
are equivalent:
\ben 
\item[i)]$(A,\cS)$ is compatible. 
\item[ii)]$(A,P_{\overline{R(A)}}(\cS))$ is compatible. 
\item[iii)]$(A,\cW)$ is compatible, for every subspace $\cW$ such that
$P_{\overline{R(A)}}(\cW)=P_{\overline{R(A)}}(\cS)$. \een
\end{cor}

\bdem 
As proved before, $\cS +N(A)$
is closed if and only if $P_{\overline{R(A)}}(\cS)$ is closed, 
so that item ii) makes sense.
\ben\item[i)$\leftrightarrow $ ii):]  As $\cS +N(A)= P_{\overline{R(A)}}(\cS)\oplus N(A)$ then
$\cS+A^{-1}(\cS^\perp)=\cS+N(A) +A^{-1}(\cS^\perp)=
P_{\overline{R(A)}}(\cS)\oplus A^{-1}(\cS^\perp)$, because
$N(A)\inc A^{-1}(\cS^\perp)$. Therefore $(A,\cS)$ is compatible if
and only if $(A,P_{\overline{R(A)}}(\cS))$ is compatible.

\item[ii)$\leftrightarrow$ iii):] Using that  i)$\leftrightarrow $ ii) for $\cS $ and $\cW$, 
$(A,P_{\overline{R(A)}}(\cW))$ is compatible \iiff 
$(A,P_{\overline{R(A)}}(\cW))$ is compatible if and only if
$(A,P_{\overline{R(A)}}(\cS))$ is compatible \iiff 
$(A,\cS)$ is compatible. 
%Observe that if
%$P_{\overline{R(A)}}(\cW)=P_{\overline{R(A)}}(\cS)$ then
%$A^{-1}(\cS\orto)=A^{-1}(\cW\orto)$ and $\cW+N(A)$ is closed so
%that $(A,\cW)$ is compatible if and only if
%$(A,P_{\overline{R(A)}}(\cW))$ is compatible; but
%$P_{\overline{R(A)}}(\cW)+A^{-1}(\cW\orto)=P_{\overline{R(A)}}(\cS)+A^{-1}(\cS\orto)$
%so \\
%$(A,P_{\overline{R(A)}}(\cW))$ is compatible if and only if
%$(A,P_{\overline{R(A)}}(\cS))$ is compatible. 
\een 
\edem

\medskip
\begin{rem}\label{4} \rm 
The pair
$(A,\cS)$ is compatible if and only if $A^{1/2}(\cS)$ is closed in
$R(A^{1/2})$ and $R(A^{1/2})=\overline{A^{1/2}(\cS)}\cap
R(A^{1/2})\oplus A^{1/2}(\cS)^\bot\cap R(A^{1/2})$.
This type of decomposition will be simplified later. For the
proof, observe first that $(A,\cS)$ is compatible if and only if
$\cH=\cS+A^{-1}(\cS\orto)$. Applying $A^{1/2}$ to both sides of this
equality we get $R(A^{1/2})=A^{1/2}(\cS)+ A^{1/2}(\cS)\orto\cap
R(A^{1/2})$. From the proposition above, $A^{1/2}(\cS)$ is closed
in $R(A^{1/2})$, so that $A^{1/2}(\cS)=\ov{A^{1/2}(\cS)}\cap
R(A^{1/2})$. The converse is similar.
\end{rem}
\bigskip

\begin{cor}\label{5.1}
The following conditions are equivalent:

i) If $\cW=A^{-1/2}(\overline{A^{1/2}(\cS)})$ then $(A,\cW)$ is
compatible.

ii) $R(A^{1/2})=\overline{A^{1/2}(\cS)}\cap R(A^{1/2})\oplus
A^{1/2}(\cS)^\bot\cap R(A^{1/2})$.

iii) There exists a solution $Q$ of $A^{1/2}X=P_\cM A^{1/2}$,
where $\cM=\overline{A^{1/2}(\cS)}$.
\end{cor}

\bdem \ \ben \item[i)] $\rightarrow$ ii) $(A,\cW)$ is compatible
if and only if $\cH=\cW+A^{-1}(\cW^\bot)$; as before, it follows
that $R(A^{1/2})=A^{1/2}(\cW)+A^{1/2}(\cW)^\bot\cap R(A^{1/2})$.

Observe that $A^{1/2}(\cW)=\overline{A^{1/2}(\cS)}\cap R(A^{1/2})$
and since 
$A^{1/2}(\cS)\inc A^{1/2}(\cW)$
$\inc\overline{A^{1/2}(\cS)}$, we get
$A^{1/2}(\cS)^\bot=A^{1/2}(\cW)^\bot$. Thus, $R(A^{1/2})=
\overline{A^{1/2}(\cS)}\cap R(A^{1/2})+A^{1/2}(\cS)^\bot\cap
R(A^{1/2})$.
\medskip

The converse is similar.
\medskip

\item[ii)]$\leftrightarrow$ iii) It is a consequence of Douglas
theorem. \een \edem

\medskip
\begin{lem}\label{CH3.9}
If $A \in \ca ^+$ then  the following conditions are equivalent:
\ben \item [i)] $R(PAP)$ is closed; \item [ii)]$A^{1/2}(\cS )$ is closed;
\item [iii)]$A(\cS )$ is closed.\een
Any of the conditions above implies that the pair $(A,\cS)$ is
compatible. In particular, if $A(\cS)$ is finite dimensional then
$(A,\cS)$ is compatible.
\end{lem}

\bdem \   Since $A\rai (\cS) = R(A\rai P)$ and $PAP = (A\rai
P)^*A\rai P$, we get the equivalence between conditions i) and ii).
Suppose that $R(PAP)$ is closed. Observe that $A(S) = R(AP)$ and
$R(AP)$ is closed \iiff $R(PA)$ or equivalently if $R(PA^2 P)$ is
closed. Note that $(PAP)^2 \le PA^2P$ and
$$N(PAP)^2 =N(PA^2P) = \cS\orto \oplus (\cS \cap N(A)).
$$
Since $ PA^2P\ge (PAP)^2 >0$ in $(N(PAP)^2)\orto$ we get that
$R(PA^2P)$ is closed. The converse is similar.
We have already proved that if $R(PAP)$ is closed then $(A,\cS)$
is compatible. \edem
\medskip
\begin{rem}\rm
The lemma shows, in particular, that in finite dimensional Hilbert
spaces,  compatibility is automatically satisfied. However, an
efficient algorithm for finding  every element of $P(A,\cS)$ is not
known.
\end{rem}

\medskip
\noi
Next, we show that $P_{A,S}$ is, modulo the orthogonal projection 
onto $\cN = N(A)\cap \cS$, a reduced solution of a Douglas-type
equation:
\medskip
\begin{pro}\label{CH4.3}
If the pair \com, let 
$\cM =\overline{A\rai (\cS) }$. Then the reduced solution $Q$ of the equation
\beq\label{PAPXPA} (PAP)X=PA
\end{equation}
coincides with the reduced solution of
\beq\label{araiP} (A\rai P)X = P_\cM A\rai .
\end{equation}
Moreover, $Q = P_{A, \cS\ominus \cN}$ and $\PA = Q + P_\cN$.
\end{pro}

\bdem Let $Q$ be the reduced solution of equation (\ref{PAPXPA}).
Observe first that $N(PAP) = N(A^{1/2}P) = \cN + \cS \orto$ and, 
therefore, $\overline{R(PAP)} = \cS
\ominus \cN$. By the definition of reduced solution,
$R(Q) \inc \overline{R(PAP)} = \cS \ominus \cN$
and $N(Q) = N(PA) =(AS) \orto$.
If $z \in \cS \ominus \cN$, then $PAPQz = PAz = PAPz$ and,
since $PAP$ is injective on $\cS \ominus \cN$, we get $Qz = z$. 
Since the pair \com,
it follows that $\cH = (\cS \ominus \cN) \oplus A^{-1}(\cS \orto)$,
so that 
$Q$ is the projection onto $R(Q) = \cS \ominus \cN$ with
$$N(Q) = A^{-1}(\cS \orto) \inc A^{-1}((\cS \ominus \cN) \orto)
= A^{-1}(R(Q) \orto).
$$
By Krein's lemma it follows that $Q \in P(A, \cS\ominus \cN )$.
Observe also that $(\cS\ominus \cN) \cap N(A) = \{ 0 \}$, so that
$P(A, \cS\ominus \cN )$ consists of a single element, namely
$P_{A, \cS\ominus \cN}$. Since $P_\cN$ is an $A-$ Hermitian
projection onto $\cN$, it follows that $Q + P_\cN =
P_{A, \cS\ominus \cN } + P_{\cN} = \PA$.

Let us prove that $Q = \PA - P_{\cN}$ is the reduced solution of
the equation (\ref{araiP}).
%For this, observe first the identity
Note  that
$(AP)Q = A^{1/2}P_{\cM}A^{1/2} =A \PA =AP \PA $. 
%, which is very easy to check (or see \cite {[Sz]}), so that
Hence $A^{1/2}(P_{\cM}A^{1/2}- A^{1/2}PQ)= \{ 0 \}$ and
$R(P_{\cM}A^{1/2}- A^{1/2}PQ) \inc N(A)$.
But, also, 
$$ R\big(P_{\cM}A^{1/2}- A^{1/2}PQ \big) \inc \overline {R(A)},
$$
so that $P_{\cM}A^{1/2}- A^{1/2}PQ = \{ 0 \},$
which says that $Q$ is a solution of (3).
In order to see that $Q$ is the reduced solution, observe
that $N(Q) = A^{-1}(\cS \orto) = N(P_{\cM}A^{1/2})$.
\edem

\medskip

If, in addition to the hypothesis of the proposition, $R(PAP)$
is supposed
to be closed, then 
$$
Q=(PAP)^\dagger PA = (A\rai P )^\dagger P_\cM A\rai =
(A\rai P)^\dagger A\rai .
$$
In fact, $PAP$ has closed range if and only if $A\rai P$ has, so that
the Moore-Penrose
inverses of these operators are bounded and, by the comments following
Douglas theorem, the reduced solution of $(PAP)X = PA$
is $(PAP)^\dagger PA$ and that of
$(A\rai P)X = P_{\cM}$ is $(A\rai P )^\dagger P_\cM A\rai$;
finally, $(A\rai P )^\dagger P_\cM  =
(A\rai P)^\dagger $ because both operators satisfy the defining
equations of the Moore-Penrose pseudoinverse of $A\rai P$.

\medskip

Concerning the minimal properties of $\PA$ in $P(A,\cS)$, 
mentioned in the Introduction,
we describe two of them. First, $\|\PA\| \le \|Q\|$ for all $Q \in P(A,\cS)$, 
but, in general,
it is not the unique element of $P(A,\cS)$ with this property (see \cite {[Sz]}). 
In order to describe the second property, we introduce some notation: 
if $T\in L(\cH,\cK)$, $\cS$ is a closed subspace of $\cH$ and $x\in \cH$,
then a \it $(T,\cS)$-interpolant \rm of $x$ is an element of the set 
spl$(T,\cS,x)=\{z\in x+\cS:\|Tz\|=\inf_{s\in\cS}
\|T(x+s)\|\}$. If $A=T^*T$, the following conditions holds:
\ben
\item[(1)] spl $(T,\cS,x)=\{z\in x+\cS:|z|_A=\inf\{|x-s|_A:s\in \cS\}\}$, 
if $|.|_A$
denotes the seminorm induced by $A$, i. e., $|x|_A ^2 = \api Ax, x \cpi$;

\item[(2)] spl $(T,\cS,x)$ is not empty for all $x\in\cH$ if and only if $(A,\cS)$ is compatible;

\item[(3)] spl $(T,\cS,x)$ has a unique element for all $x\in \cH$ if and only if $(A,\cS)$ is compatible
and $\cS\cap N(A)=0$;

\item[(4)] If $(A,\cS)$ is compatible and $x\in \cH$ then $(1-P_{A,\cS})x$ is the unique element of
spl $(T,\cS,x)$ with minimal norm. 
\een
\noi
The proofs of these facts can be found in \cite{[CMS3]}.

\section{Operator ranges}

In this section we recall a well-known construction of a Hilbertian structure
on the range of an operator (see \cite{[Br]}, \cite{[FW]} or \cite{[An]}).
We include some new facts which will be prove useful in the following section.

\begin{num}\rm
As we have already seen any $A\in L(\cH)^+$ induces a semi-inner
product on $\cH$, by means of $(x,y)_A=\langle Ax,y\rangle$, 
$x,y\in\cH$. 
Denote by $\cK = \ov{R(A)}$ and $P_\cK \in \ca$ the orthogonal 
projection onto $\cK$. 
Since $\cK = N(A) \orto$, we can 
define on $\cK$ the inner product 
$$(x, y)_A = \langle Ax,y\rangle~,\quad x,y\in\cK.
$$
This inner product induces the norm $|x|_A = 
(x, x)_A\rai  = \langle Ax,x\rangle \rai$, $x\in\cK$.
Let $\cH_A$ be the completion of the inner product space 
$(\cK,(~,~)_A)$. We assume that $\cK \inc \cH_A$. 

%
%Consider the quotient space $\cH/N(A)$ and the natural
%map $\pi:\cH\to\cH/N(A)$, defined by $\pi(x)=x+N(A)=\ov x$, the
%quotient class of $x\in\cH$. Observe that $\cH/N(A)$ is isomorphic
%to $\ov{R(A)}$. Define on $\cH/N(A)$ the inner product
%$(\ov{x},\ov{y})_A=\langle Ax,y\rangle~,\quad x,y\in\cH$. Then
%$(\ov{x},\ov{y})_A$ is well defined and induces a norm
%$|\ov{x}|_A=(\ov{x},\ov{x})_{A}= \|A^{1/2}x\|^2$. Let $\cH_A$ be the
%completion of the inner product space $(\cH/N(A),(~,~)_A)$; or
%equivalently, the completion of $\left(\ov{R(A)}, (~,~)_A\right)$.
%Define $\varphi=i\circ\pi$, where
%$$
%i: (\cH/N(A),(~,~)_A)\to (\cH_A, (~,~)_A)
%$$
%is the natural inclusion. 
Then, the projection $P_\cK$ induces a
map $\varphi :\cH\to\cH_A$, defined by $
\varphi (x) = P_\cK x $, $x\in\cH$. Note that  $\varphi$ has  
nullspace $\cK \orto = N(A)$. 

For every subspace $\cS$ of $\cH$, $\varphi(\cS)$ is a 
subspace in $\cH_A$. By $\ov{\varphi(\cS)}^{\cH_A}$ 
we denote its closure in $(\cH_A,(~,~)_A)$. 
%(where $\ov{\varphi(\cS)}^{\cH_A}$ is the closure in $(\cH_A,(~,~)_A)$)
Then the orthogonal projection  from $\cH_A$ onto
$\ov{\varphi(\cS)}^{\cH_A}$ always exists,  even if the 
original pair $(A,\cS)$ is not compatible.

The relative position in $\cH$ between $\cS$ and $N(A)$ obviously
affects the ``size" of the projection. In this section we
deduce conditions on this projection in order to obtain the compatibility of
$(A,\cS)$.
\end{num}

\begin{num}\rm
The construction of $\cH_A$ can be performed in the context of
operator ranges. We refer the reader to the papers by Dixmier
\cite{[Di]} and Fillmore and Williams \cite{[FW]} and to Ando's
book \cite{[An]}.

Consider  $T\in L(\cH,\cK)$. The range of $T$ can be given a Hilbert space
structure $(R(T),\langle~,~\rangle_T)$ in a unique way, such that
$T$ becomes a coisometry from $(\cH,\langle~,~\rangle)$ to
$(R(T),\langle~,~\rangle_T)$ (see \cite{[An]}).
More precisely, as $T:N(T)\orto\to R(T)$ is a bijection, define
$\langle Tx,Ty\rangle_T=\left\langle P_{N(T)\orto}x,
P_{N(T)\orto}y\right\rangle$, for $x,y\in\cH$.  For $u\in R(T)$,
denote $\|u\|_T=\langle u,u\rangle^{1/2}_T$. The key fact is that
the operator $T:(\cH,\langle~,~\rangle)\to
(R(T),\langle~,~\rangle_T)$ is a coisometry. 
%The norm $\|u\|_T=\langle u,u\rangle_T^{1/2}$ verifies 
Observe that $\|u\|_T=\min \{\|a\|:Ta=u\}, \hbox{ for all } u\in R(T)$, because
$Ta=TP_{N(T)\orto} a$ and $\|P_{N(T)\orto}a\|\le \|a\|$.
Therefore, $T:(\cH,\langle~,~\rangle )\to (R(T),\langle
~,~\rangle_T)$ is bounded, i.e., $\|Ta\| _T\le \|a\|$, for all
$a\in \cH$; also, for each $u\in R(T)$, there is a unique $a\in
N(T)\orto$ such that $Ta=u \quad\hbox{and}\quad \|a\|=\|u\|_T$.
As in the introduction, we use the notation
$B(T)=(R(T),\langle~,~\rangle_T)$.
In \cite{[Di]} and \cite{[FW]} a number of characterizations of
operator ranges are given. One of them, establishes that a
subspace $\cR$ of $\cH$ is the range of a bounded operator if and
only if there is an inner product $\langle ~,~\rangle '$ on $\cR$
such that $(\cR,\langle ~,~ \rangle ')$ is a Hilbert space and
$\|x\|' \ge \|x\|$ for all $x\in \cR$. (see \cite{[FW]}, Theorem
1.1). More precisely, given $T\in L(\cH,\cK)$, consider
$T_1=(T|_{N(T)\orto})^{-1}$, $T_1:R(T)\to N(T)\orto$, and define
$
\langle u, v\rangle '=\langle u, v\rangle +\langle
T_1u,T_1v\rangle , \hbox{~for~} u,v\in R(T)~.
$
Then $\langle~,~\rangle '$ is complete and $\|u\|'\ge \|u\|$ for
all $u\in \cR$. In fact, the inner products $\langle ~,~\rangle '$
and $\langle ~,~\rangle_T$ are equivalent: first, observe that if
$Tx=u$ for $x\in\cH$ then
$T_1u=T_1Tx=T_1TP_{N(T)\orto}x=P_{N(T)\orto}x$. Therefore,
$\|u\|'{^2}=\langle u, u\rangle '=\langle Tx,Tx\rangle '\ge
\langle T_1Tx, T_1Tx\rangle = \|P_{N(T)\orto}x\|^2=\|u\|^2_T$, so
that $\|u\|_T\le \|u\|'$. Conversely,
$\|u\|'{^2}=\|TP_{N(T)\orto}x\|^2+\|u\|^2_T\le
(\|T\|^2+1)\|u\|^2_T$. Observe also that $\langle u,v \rangle
'=\langle u,v\rangle +\langle u,v\rangle_T$, for $u,v\in \cR$.
\end{num}

\bigskip
\noi
We shall consider the construction above for any positive operator
on $\cH$. More precisely, given $A\in L(\cH)^+$ consider
$R(A^{1/2})$ with the norm induced by $\langle~,~\rangle_
{A^{1/2}}$, i.e. the space $B(A^{1/2})$. The next lemma shows that
$A$ provides an isometric isomorphism between $\cH_A$ and
$B(A^{1/2})$. It should be mentioned that the subtle
relationship between $R(A)$ and $R(A^{1/2})$ is fundamental in
this and all remaining results.
\medskip
\begin{lem}\label{15}
Given $A\in L(\cH)^+$
$$
A|_{\ov{R(A)}}:\left(\ov{R(A)},(~,~)_A\right)\to B(A^{1/2})
$$
is an isometry with dense image and then it admits a unitary
extension
$$
A': \left(\cH_A,(~,~)_A\right)\to B(A^{1/2}).
$$
\end{lem}

\bdem Denote, as before, $\cK = \ov{R(A)}$. For all $x\in \cK$, it holds
$$\|Ax\|_{A^{1/2}}=\|P_{\cK }A^{1/2}x\|=\|A^{1/2}x\|= |x|_A.
$$
Also, if $x\in \cK $ and $\|Ax\|_{A^{1/2}}=0$ then
$|x|_A=\|A^{1/2}x\|=0$ so that $x\in N(A^{1/2})=N(A)$, and $x=0$.
It remains to prove that the image of $A|_{\cK }$
is dense in $B(A^{1/2})$: since $R(A^{1/2})$ is dense in 
$\cK $, for any $x \in \cH$
there exists a sequence $\{x_n\}$ in $\cH$ such that
$A^{1/2}x_n\to P_{\cK }x$, which means that $Ax_n\to
u$ in $B(A^{1/2})$, for any $u \in R(A^{1/2})$. 
Then $A|_{\cK } :\cK \to R(A)\inc
B(A^{1/2})$ admits a unitary extension from the completion of
$\cK $, namely $\cH_A$, onto $B(A^{1/2})$. \edem

\medskip
\begin{rem}\rm
More generally, for $t\in[0,1]$ consider $A^t$ and define in $R(A^t)$ 
the inner product
$\langle A^tx,A^ty\rangle_{A^t}=\langle
P_{\cK }x,P_{\cK }y\rangle$ . Observe
that $\ov{R(A^t)}=\cK $, for all $t\in [0,1]$. Denote
$\langle A^tx,A^ty\rangle_t=\langle A^tx, A^ty\rangle_{A^t}$, for
$x,y\in\cH$, $\|A^tx\|_t=\|A^tx\|_{A^t}=\|A^{t/2}x\|$,
$(x,y)_t=(x,y)_{A^t}$ and $\cH_t=\cH_{A^t}$. Then
$|x|_t=|x|_{A^t}= \|A^{t/2}x\|$.
\end{rem}
\medskip
\noi
As before, we get:

\medskip
\begin{cor} Given $A\in L(\cH)^+$, the operator
$$
A^t|_{\ov{R(A)}}:(\ov{R(A)},(~,~)_t)\to B(A^{t/2})
$$
is an isometry with dense image and it admits a unitary extension
$$
(A^t)': (\cH_t,(~,~)_t)\to B(A^{t/2})~.
$$
\end{cor}

\bdem Straightforward. \edem

\noi
We have the following commutative diagram:
$$
\begin{array}{ccccccc}
\cH & \stackrel{\varphi}\longrightarrow & (\ov{R(A)},(~,~)_A)
&\stackrel{i}\longrightarrow & (\cH_A,(~,~)_A )\\
& &\Big\downarrow\vcenter{\rlap{$A$}} &  & \Big\downarrow\vcenter{\rlap{$A'$}}\\
\cH & \stackrel{A}\longrightarrow & (R(A),\langle
~,~\rangle_{A^{1/2}}) &\stackrel{j}\longrightarrow & B(A^{1/2})
\end{array}
$$
which relates two maps from $\cH$ into two Hilbert spaces
associated to the operator $A$. By the lemma above these two spaces are
isometrically isomorphic. Observe that the images of the subspace
$\cS$ of $\cH$ in $\cH_A$ and $B(A^{1/2})$ are, respectively,
$\phi(A)$ and $A(\cS)$.

The next result, which has been proved by Barnes for injective
operators (see \cite{[Bar]}), characterizes the operators $B\in
L(\cH)$ which can be extended to $L(B(A^{1/2}))$.

\bigskip

\begin{lem}\label{16}
Consider $B\in L(\cH)$. There exists  $\widetilde B\in
L(B(A^{1/2}))$ such that $\wtil B A=AB$ if and only if
$B(N(A))\inc N(A)$ and $R(B^*A^{1/2})\inc R(A^{1/2})$.
In this case, such an operator is unique.
\end{lem}

\bdem Let $\wtil B\in L(B(A^{1/2}))$ such that $\wtil BA=AB$; if
$x\in N(A)$ then $ABx=0$ so that $Bx\in N(A)$ and $B(N(A))\inc
N(A)$. Since $\wtil B\in L(B(A^{1/2}))$, there exists $C>0$ such
that $\|\wtil BAx\|_{A^{1/2}}\le C\|Ax\|_{A^{1/2}}$ for all $x\in
\cH$; equivalently, $ \|ABx\|_{A^{1/2}}\le C\|Ax\|_{A^{1/2}}$. By
definition  of $\|~\|_{A^{1/2}}$, this means
$$
\|P_{\ov{R(A)}}A^{1/2}Bx\|\le C\|P_{\ov{R(A)}}A^{1/2}x\|
\quad \hbox{  or  } \quad 
\|A^{1/2}Bx\|\le C\|A^{1/2}x\| , 
$$ 
because $R(A^{1/2})\inc
\ov{R(A)}$. By Douglas theorem, the last inequality is equivalent to
$R(B^*A^{1/2})\inc R(A^{1/2})$. Conversely, if these conditions
hold, it is easy to see that $\wtil B$ can be defined in $R(A)$
and extended to a bounded operator in $B(A^{1/2})$. If there
exists $C\in L(B(A^{1/2}))$ such that $CA=\wtil BA$ then $C$ and
$\wtil B$ coincide in $R(A)$, which is dense in $B(A^{1/2})$, so
that $C=\wtil B$. \edem
\bigskip

\def\QA{Q_{A, \cS}}
\begin{num}\rm
Given a subspace $\cW$ of $B(A^{1/2})$ the closure (resp. the
orthogonal complement) of $\cW$ in $B(A^{1/2})$ is denoted
$\ov{\cW}'$ (resp. $\cW^{\bot'}$). If $\cS$
is a closed subspace of $\cH$ and $\cW=A(\cS)$, then
$\cM = \ov{\cW}'=A^{1/2}\left(\ov{A^{1/2}(\cS)}\right)$ and
$$
\cM^{\bot'}=\cW^{\bot'}=A^{1/2}(A^{1/2}(\cS)\orto)=\cS\orto\cap R(A^{1/2}).
$$
From now $\QA \in L(B(A^{1/2}))$ denotes the 
orthogonal projection onto $\cM$. Then
$R(\QA )=\cM=A^{1/2}(\ov{A^{1/2}(\cS)})$ and $N(\QA )=\cM^{\bot'}=
\cS\orto\cap R(A^{1/2})$. Observe that $A(\cS)\inc R(\QA )$ and
$A(A^{-1}(\cS\orto))=\cS \orto \cap R(A)\inc N(\QA )$.
\end{num}

\medskip
\begin{lem}\label{17}
It holds $\ov{\cS\orto\cap R(A)}'=\cM^{\bot'}$ if and only if
$A(\cS)+\cS\orto \cap R(A)$ is dense in $B(A^{1/2})$.
\end{lem}

\bdem $\ov{A(\cS)+\cS\orto\cap R(A)}'=\cM +\ov{\cS\orto\cap
R(A)}'$ because $\cS\orto\cap R(A)\inc \cM^{\bot '}$. Then
$\ov{\cS\orto\cap R(A)}'=\cM ^{\bot '}$ if and only if
$A(\cS)+\cS\orto\cap R(A)$ is dense in $B(A^{1/2})$. \edem
\bigskip

\section{Compatibility and operator ranges}

We have now the tools for proving the relationship between the
compatibility of $A$ with $\cS$ and the properties of the
orthogonal projection $\QA \in L(B(A\rai ))$ onto $\cM=\ov{A(\cS)}'$. We start with a
technical result.

\medskip
\begin{pro}\label{18}
Given $A\in L(\cH)^+$ the following conditions are equivalent:
\ben \item[i)] $(A,\cS)$ is compatible.
\item[ii)]
$A^{1/2}(\cS)+A^{1/2}(\cS)\orto\cap R(A^{1/2})$ is closed in
$R(A^{1/2})$ and $A^{1/2}(\cS)\orto\cap R(A^{1/2})$ is dense in
$A^{1/2}(\cS)\orto\cap \ov{R(A^{1/2})}$.
\item[iii)] $R(A)=\cM\cap
R(A)+\cM^{\bot'}\cap R(A)$ and $A(\cS)$ is closed in $R(A)$ in the
topology of $B(A^{1/2})$. \een
\end{pro}

\bdem 
Denote $\cT = A^{1/2}(\cS)\orto\cap R(A^{1/2})$. 
\ben \item[i)]$\rightarrow$ ii) By Remark \ref{4}
$(A,\cS)$ is compatible if and only if
$R(A^{1/2})=A^{1/2}(\cS)+\cT$, so
that $A^{1/2}(\cS)+\cT$ is closed in
$R(A^{1/2})$. Also
$$\ov{R(A^{1/2})}=\ov{A^{1/2}(\cS)}+\ov{\cT}\inc
\ov{A^{1/2}(\cS)}+A^{1/2}(\cS)\orto\cap\ov{R(A^{1/2})}\inc\ov{R(A^{1/2})} .
$$ 
Therefore $\cT$ is dense in
$A^{1/2}(\cS)\orto\cap \ov{R(A^{1/2})}$.

\item[ii)]$\rightarrow$ iii) By assumption, it holds 
$$
\barr{rl}
A^{1/2}(\cS )+ \cT & = \ov{A^{1/2}(\cS)+\cT} \cap R(A^{1/2})\\ 
& = \left( \ov{A^{1/2}(\cS)}+ \ov { \cT } \right)\cap R(A^{1/2}) \\
& = \ov{R(A^{1/2})}\cap R(A^{1/2})=R(A^{1/2}),
\earr
$$ 
because $\ov{ \cT }=A^{1/2}(\cS)\orto\cap\ov{R(A^{1/2})}$.
Then $R(A^{1/2})=A^{1/2}(\cS)+\cT $
and 
$$
\barr{rl}
R(A)& =A^{1/2}(R(A^{1/2}))= A(\cS)+A^{1/2}(\cT )\\
& \inc A(\cS)+\cM{\orto}'\cap R(A) \\
& \inc \cM \cap R(A)+\cM{\orto}'\cap R(A)\inc R(A).
\earr
$$
Therefore $R(A)=\cM\cap R(A)+\cM\orto{'}\cap R(A)$ and $\cM\cap
R(A)=A(\cS)$.

\item[iii)]$\rightarrow$ i) Straightforward.

\een \edem

\medskip
\begin{rem}\label{19} \rm
Condition ii) is equivalent to ii') $A(\cS)+\cS\orto\cap R(A)$ is
closed in $R(A)$ under the topology of $B(A^{1/2})$ and
$\ov{\cS\orto\cap R(A)}'=\cM^{\bot'}$. In fact, from the proof of
ii)$\to$ iii) we get that $R(A)=A(\cS)+\cM{\orto}'\cap
R(A)=A(\cS)+\cS\orto\cap R(A)$ so that $A(\cS)+\cS\orto\cap R(A)$
is closed in $(R(A),\langle~,~\rangle_{A^{1/2}})$ and
$\ov{\cS\orto \cap R(A)}'=\cM^{\bot '}$. The converse is similar.
\hfill$\square$
\end{rem}
\bigskip

\noi
In the last part of the paper, we relate the compatibility of the pair $(A,\cS)$
with the existence of certain projections in $B(A^{1/2})$. As
before, $\cM=\ov{A(\cS)}'=A^{1/2}(\ov{A^{1/2}(\cS)})$.

\medskip
\begin{teo}\label{20}
If $(A,\cS)$ is compatible then there exists $\wtil{P}_{A,\cS}\in
L(B(A^{1/2}))$ such that $\wtil{P}_{A,\cS}A=AP_{A,\cS}$. Moreover,
$\wtil{P}_{A,\cS}=\QA $.
\end{teo}

\bdem If $(A,\cS)$ is compatible then, by proposition 3.10,
$R=P_{A,\cS\ominus\cN}$ is the reduced solution of
$$
A^{1/2}P X=P_{\ov{A^{1/2}(\cS)}}A^{1/2}
$$
and $P_{A,\cS}=R+P_\cN$ where $\cN=\cS\cap N(A)$. Observe that
$A^{1/2}P_{A,\cS} = A^{1/2}P P_{A,\cS}=A^{1/2}
P R =P_{\ov{A^{1/2}(\cS)}}A^{1/2}$  because $R(P_{A,\cS})=\cS$. Therefore
$P_{A,\cS}$ verifies
$A^{1/2}P_{A,\cS}=P_{\ov{A^{1/2}(\cS)}}A^{1/2}$ so that
$P^*_{A,\cS}A^{1/2}=A^{1/2}P_{\ov{A^{1/2}(\cS)}}$ and then
$R(P^*_{A,\cS}A^{1/2})\inc R(A^{1/2})$. In order to apply lemma
\ref{16}, let $x\in N(A)$ and observe that
$$A^{1/2}P_{A,\cS}x=P_{\ov{A^{1/2}(\cS)}}A^{1/2}x=0,
$$
because
$N(A)=N(A^{1/2})$, so $P_{A,\cS}x\in N(A)$. By
\ref{16}, there exists $\wtil{P}_{A,\cS}\in L(B(A^{1/2}))$ such
that $\wtil{P}_{A,\cS}A=AP_{A,\cS}$; now,
$\wtil{P}_{A,\cS}(R(A))=A(R(P_{A,\cS}))=A(\cS)$, so that
$A(\cS)\inc R(\wtil{P}_{A,\cS})$ and $\cM=\ov{A(\cS)}'\inc
R(\wtil{P}_{A,\cS})$. Also, $A(N(P_{A,\cS}))\inc
N(\wtil{P}_{A,\cS})\inc \cM^{\bot '}$ which implies 
$\ov{\cS\orto\cap R(A)}'\inc N(\wtil{P}_{A,\cS})\inc\cM^{\bot'}$.
But $\ov{\cS\orto\cap R(A)}'=\cM^{\bot'}$ because $(A,\cS)$ is
compatible and the proposition and remark above apply. Then
$R(\wtil{P}_{A,\cS} )=\cM$, so that $\wtil{P}_{A,\cS} = \QA$. \edem

\noi
The next theorem gives a simple characterization of compatibility:
\smallskip
\begin{teo}\label{22}
$(A,\cS)$ is compatible if and only if $\QA (R(A))=A(\cS)$.
\end{teo}
\bdem If $(A,\cS)$ is compatible then, by theorem \ref{20},
$\wtil{P}_{A,\cS}=\QA $, so that 
$$
\QA (R(A))=\QA A(\cH)=AP_{A,\cS}(\cH) =A(\cS).
$$
Conversely, if $\QA (R(A))=A(\cS)$, any $x\in R(A)$
decomposes as $x=x_1+(I-\QA )x$, where $x_1=As$ for some $s\in\cS$;
then $(I-\QA )x\in N(\QA )\cap R(A)=\cS\orto\cap R(A)$ and
$R(A)=A(\cS)+\cS\orto\cap R(A)$. Then, $\cH=\cS+A^{-1}(\cS\orto)$,
which shows that $(A,\cS)$ is compatible. \edem

\noi
Denote $A^\sh=(A|_{\ov{R(A)}})^{-1}:R(A)\to\ov{R(A)}$.
\medskip
\begin{lem}\label{21}
The projection $\QA $ satisfies $\QA (R(A))\inc R(A)$ if and only if $R(A)=\cM\cap
R(A)+\cM^{\bot'}\cap R(A)$; in this case $\QA (R(A))= \cM\cap R(A)$.
Moreover, $A^\sh \QA A:\cH\to \cH$ is a bounded projection if and
only if $\cM\cap R(A)$ is closed in $R(A)$ (under the topology of
$\cH$).
\end{lem}

\bdem Observe that, by the definition of $\cM$, $\QA (R(A))\inc R(A)$
if and only if $\QA x\in \cM\cap R(A)$, for all $x\in R(A)$. Then
$(I-\QA )x\in \cM^{\bot'}\cap R(A)$ so that $R(A)=\cM\cap
R(A)+\cM^{\bot'} \cap R(A)$. On the other hand, it always holds
$\cM\cap R(A)\inc \QA (R(A))$. Then $\QA (R(A)) =\cM\cap R(A)$. The
converse is similar. If $\QA (R(A))\inc R(A)$ then $A^\sh
\QA A:\cH\to\cH$ is well defined and it is obviously a projection.
Let us prove that it is bounded. For this, observe that $N(A^\sh
\QA A)=N(\QA A)=A^{-1}(\cS\orto)$ is closed and, also, $R(A^\sh
\QA A)=A^\sh \QA (R(A))=A^\sh(\cM\cap R(A))$ is closed, because $\cM\cap
R(A)$ is closed in $R(A)$. This proves that $A^\sh \QA A$ is bounded.
\edem

\bigskip
\noi
Consider now the following subalgebra of $L(\cH)$:
$$
L(\cH)^A=\{T\in L(\cH):T(N(A))\inc N(A)~~\hbox{ and }~~
R(T^*A^{1/2})\inc R(A^{1/2})\}.
$$
By lemma \ref{16} the elements of $L(\cH)^A$ induce operators on
$B(A^{1/2})$ by means of %the map

$$
\theta : L(\cH)^A \to L(B(A^{1/2}))\quad \hbox{ given by } \quad
T \mapsto \theta(T)=\wtil{T} , \ \ T \in \ca^A , 
$$
where $\wtil{T}Ax=ATx$, for all $x\in\cH$.

\medskip
\begin{teo}\label{23}
Given a closed subspace $\cW$ of $B(A^{1/2})$ and $Q_\cW\in
L(B(A^{1/2}))$ the orthogonal projection onto $\cW$, then
$\theta^{-1}(\{Q_\cW\})$ is non empty if and only if $(A,\cS)$ is
compatible, where $\cS$ is any closed subspace of $\cH$ such that
${A(\cS)}$ is dense in  $\cW$.
\end{teo}

\bdem If there exists $\cS$ such that $(A,\cS)$ is compatible and
${A(\cS)}$ is dense in $\cW$ then, by theorem \ref{20}, there
exists $\wtil{P}_{A,\cS}\in L(B(A^{1/2}))$ such that
$\wtil{P}_{A,\cS} =Q_\cW$. Therefore $\theta(P_{A,\cS})=Q_\cW$.

Conversely, if $\theta^{-1}(Q_\cW)$ is non empty, then there exists
$T\in L(\cH)^A$ such that $\wtil T=Q_\cW$ and $\wtil TA=AT$; then
$Q_\cW(R(A))\inc R(A)$. By lemma \ref{21} this inclusion is
equivalent to $R(A)=\cW\cap R(A)+\cW^{\bot'}\cap R(A)$ and in this
case $Q_\cW(R(A)) =\cW\cap R(A)$. Then $\ov{\cW\cap R(A)}'=\cW$
because $R(A)$ is dense in $B(A^{1/2})$.
Again by lemma \ref{21}, $\cW\cap R(A)$ is closed in $R(A)$
because $A^\sh Q_\cW A=P_{\ov{R(A)}}T$ is a bounded projection in
$\cH$. Then $\cS=A^{-1}(\cW\cap R(A))$ is a closed subspace of
$\cH$ such that $A(\cS)=\cW\cap R(A)$ in $\cH$ because $\cW\cap
R(A)$ is closed, so that $\ov{A(\cS)}'=\cW$.
Applying theorem \ref{22} to $Q_\cW $, since
$Q_\cW(R(A))=A(\cS)$, with $\ov{A(\cS)}'=\cW$ we obtain that
$(A,\cS)$ is compatible. \edem

\medskip
\begin{pro}\label{24}\rm
Let $\cS$ is a closed subspace of $\cH$ and $\cW=\ov{A(\cS)}'$. If
$\theta^{-1}(\{Q_\cW\})$ is non empty, then $P(A,\cS)\inc
\theta^{-1}(\{Q_\cW\})$. Moreover $P(A,\cT)\inc
\theta^{-1}(\{Q_\cW\})$ for all closed subspaces $\cT$ of $\cH$,
such that $A(\cT)=A(\cS)$.
\end{pro}

\bdem If $R\in P(A,\cS)$ then $R=P_{A,\cS}+T$, where 
$T\in L(\cS\orto,\cN)$ (see \cite{[Sz]}). Then
$\theta(R)=\theta(P_{A,\cS})=\wtil P_{A,\cS}$ because $\wtil T
A=AT=0$.

If $A(\cT)=A(\cS)$ then $P_{\ov{R(A)}}(\cT)=P_{\ov{R(A)}}(\cS)$.
Observe that, by proposition \ref{23}, $(A,\cS)$ is compatible
because $\theta^{-1}(\{Q_\cW\})$ is non empty; therefore $(A,\cT)$
is compatible by corollary \ref{R3.6}. But, by proposition
\ref{20}, $\wtil P_{A,\cT}$ is the orthogonal projection onto
$\ov{A(\cT)}'=\cW$ so that $\theta (P_{A,\cT})=\theta
(P_{A,\cS})$, and $\theta(P(A,\cT))=\theta(P_{A,\cT})$. \edem

\medskip
\begin{rem}\label{25}\rm
If one decides to avoid the use of operator ranges with their
natural Hilbertian structure, then by remark 3.6,  
$(A,\cS)$ is compatible if and only if $A^{1/2}(\cS)$ is a closed 
subspace of $R(A^{1/2})$ which admits an orthogonal 
complement in  $R(A^{1/2})$; observe that, as a subspace of $\cH$,
$R(A^{1/2})$ is an incomplete inner product space,  unless $R(A)$ is closed. 
Therefore, the
compatibility problem is equivalent to find 
in an inner-product space $(\cD,\langle~,~\rangle)$, all 
closed subspaces of $\cD$ that admit an
orthogonal complement in $\cD$. These subspaces are called Chebyshev
subspaces in the theory of best approximation (see \cite {[De1]} for an nice treatment of
Chebyshev sets in inner product spaces).
\end{rem}

\bigskip

\noindent{IAM-CONICET and 
Departamento de Matem\'atica, FI-UBA, \\
Paseo Col\'on 850,  Buenos
Aires (1063), Argentina}

\noindent{e-mail: gcorach@fi.uba.ar}

\medskip

\noindent{IAM-CONICET and  Instituto de Ciencias, UNGS,
Los Polvorines, Argentina}

\noindent{e-mail : amaestri@ungs.edu.ar}
\medskip

\noindent{IAM-CONICET and Departamento de Matem\'atica, FCE-UNLP, \\
La Plata,  Argentina } 

\noindent{e-mail:
demetrio@mate.unlp.edu.ar}

\end{document}